\documentclass[pdflatex,sn-mathphys-num]{sn-jnl}

\usepackage{graphicx}%
\usepackage{multirow}%
\usepackage{amsmath,amssymb,amsfonts}%
\usepackage{amsthm}%
\usepackage{mathrsfs}%
\usepackage[title]{appendix}%
\usepackage{xcolor}%
\usepackage{textcomp}%
\usepackage{manyfoot}%
\usepackage{booktabs}%
\usepackage{algorithm}%
\usepackage{algorithmicx}%
\usepackage{algpseudocode}%
\usepackage{listings}%
\usepackage{cuted}%

\def\rank #1{\mathrm{rank}\,#1}

\def\bsy{\boldsymbol}

\def\wt{\widetilde}
\def\ch #1{\mathrm{Char}\,#1}
\def\w{\wedge}
\def\sfl{static feedback linearizable}
\def\dfl{dynamic feedback linearizable}
\def\stf{static feedback}
\def\df{dynamic feedback}

\def\bsy{\boldsymbol}

\def\wt{\widetilde}
\def\ch #1{\mathrm{Char}\,#1}
\def\w{\wedge}
\def\mrm #1{\mathrm{#1}}

\def\CV{{\mathcal{V}}}
\def\CH{{\mathcal{H}}}

\def\CB{{\mathcal{B}}}

\def\CH{{\mathcal{H}}}

\def\CL{{\mathcal{L}}}

\def\mcal #1{\mathcal{#1}}

\def\wt{\widetilde}

\def\b{\beta}
\def\D{\Delta}

\def\g{\gamma}
\def\G{\Gamma}

\def\k{\kappa}

\def\vf{\varphi}

\def\o{\omega}

\def\r{\rho}

\def\t{\tau}
\def\th{\theta}

\def\i{\text{\rm i}}
\def\B#1{\mathbb{#1}}

\def\P #1{\partial_{#1}}
\def\ds{\displaystyle}
\def\i{\text{\rm i}}
\def\B#1{\mathbb{#1}}

\def\P #1{\partial_{#1}}
\def\ds{\displaystyle}

\newtheorem{theorem}{Theorem}
 
\newtheorem{corollary}[theorem]{Corollary}
\newtheorem{definition}{Definition}%

\begin{document}

\title[Article Title]{Flatness of the Generic Under-Actuated 6-Degrees-of-Freedom Unmanned Surface Vessel}


\author[1]{\fnm{T.J.} \sur{Klotz}}\email{taylor.klotz@colorado.edu}

\author[2]{\fnm{P.J.} \sur{Vassiliou}}\email{peter.vassiliou@anu.edu.au}
\equalcont{These authors contributed equally to this work.}


\affil[1]{\orgdiv{Department of Applied Mathematics}, \orgname{University of Colorado}, \orgaddress{\street{} \city{Boulder}, \postcode{100190}, \state{CO}, \country{USA}}}

\affil[2]{\orgdiv{Mathematical Sciences Institute}, \orgname{Australian National University}, \orgaddress{\street{} \city{Canberra}, \postcode{2601}, \state{ACT}, \country{Australia}}}



\abstract{The flatness of the under-actuated unmanned surface vessel  with  circular hull-shape (``the hovercraft")  has long been known. We resolve the long-standing open problem: whether the generic under-actuated 6-degrees-of-freedom (dof) unmanned surface vessel (gUSV) admits a feedback linearization. Lie symmetries of the gUSV are used to construct  a dynamic feedback linearization, explicit flat outputs and flat parametrization, for all hydrodynamic parameter values. In this generic case the dynamic compensator has order 6, in contrast to the hovercraft vessel whose dynamic compensator has order 4.}




\maketitle

\section{Introduction}\label{sec1}

Feedback linearization is amoung the most consequential tools in nonlinear control: when a system can be brought, by state feedback and coordinate change\footnote{In Theorem 7.12 of \cite{CKV24} a more general diffeomorphism than a static feedback transformation is central in generating the dynamic compensator for a dynamic feedback linearization from symmetry.}, into a linear controllable form, the entire machinery of linear control -- pole placement, tracking, trajectory generation via flatness -- becomes available in closed form.  For underactuated marine vehicles,  unmanned surface vessels (USV) in particular, this question has practical significance: autonomous surface vessels are increasingly required to track prescribed trajectories with performance guarantees, and a feedback linearizing, flat parametrization would reduce planning and tracking to linear design. Yet due to its relatively complex dynamics, whether the generic 6-degrees-of-freedom (dof) unmammned surface vessel (gUSV) admits such a linearization has remained an open question.   Indeed, it is fokelore that 
the ``{\it [generic unmanned] surface vessel model} $\ldots$ {\it is notoriously not flat}'' \cite{fliessUSV}.

The gUSV studied here is the standard Fossen-type horizontal plane model with 6 states $(x,y,\theta, u,v,r)$ -- earth-frame position and heading, and body-frame surge ($u$), sway ($v$) and yaw ($r$) velocities -   with two inputs: surge force $u_1$, and yaw moment $u_2$ \cite{fossen}. This system, or special cases of it, have been the subject of multiple studies, for instance, \cite{PettNij}, \cite{LagHarmChit}, \cite{FZLS}, \cite{SLM}, \cite{fliessUSV}. In the case of the latter three references differential flatness (see Subsection \ref{flatnessSubSect}) plays a prominent role. 

The gUSV control system in question is modelled by

\begin{equation}\label{surfaceKinematics}
\begin{aligned}
&\dot{x}=u\cos\theta-v\sin\theta,\cr
&\dot{y}=u\sin\theta+v\cos\theta,\cr
&\dot{\theta}=r,\cr
&\dot{u}=\frac{1}{\gamma_1} vr-\beta_1 u+u_1,\cr
&\dot{v}=-\gamma_1 ur-\beta_2 v,\cr 
&\dot{r}= \gamma_2 uv-\beta_3 r+u_2.
\end{aligned}
\end{equation}
Here the parameters $\gamma_1, \gamma_2$ are related to the geometry of the vessel and its added mass tensor, while $\beta_1, \beta_2, \beta_3$ are positive constants quantifying hydrodynamic drag. The parameter value $\gamma_2=0$ corresponds to a circular hull--shape, sometimes referred to as ``the hovercraft'' \cite{SiraRamirez}. The hovercraft vessel is known to be flat with flat outputs $[x, y]$ 
\cite{SiraRamirez}, \cite{fliessUSV}.  References \cite{fliessUSV} and \cite{FZLS} adapt the flatness of the hovercraft to construct appoximate trajectory generation and control strategies for the generic USV in the case
$\gamma_2\neq 0$. It is therefore of interest to settle the question of the flatness of the  gUSV in the case 
$\gamma_2\neq 0$.  This is the main subject of the present paper.  We study system \eqref{surfaceKinematics} with generic values of all five parameters $\gamma_1, \gamma_2,\beta_1,\beta_2,\beta_3$, in which $\gamma_1>0$ for physically realistic vessels. For this generic family of system parameters we derive an explicit dynamic feedback linearization together with the corresponding flat outputs, which turn out to depend upon both the configuration coordinates $x, y, \th$, as well as the velocities, $u, v$ and $r$. 

\subsection{Differential flatness}\label{flatnessSubSect}
Differential flatness  is a highly sought-after property of a nonlinear control  system -- see, for instance, \cite{BeaverMalikopoulos}, \cite{rigatosBook2015}, \cite{Levine2009}.   A control system 
\begin{equation}\label{generalControlSystem}
\frac{dx}{dt}=f(t,x,u),
\end{equation}
 is said to be {\it flat} if there is a vector-valued function
$y(t)$ of the states $x(t)$, inputs $u(t)$ and finitely many derivatives of $u(t)$, which can then be inverted to express all the dynamical variables in terms of $y$ and its derivatives:
\begin{equation}\label{explicitSOLformula}
\begin{aligned}
&x(t)=\xi(y(t), \dot{y}(t),\ldots, y^{(p)}),\cr
&u(t)=\eta(y(t), \dot{y}(t),\ldots, y^{(q)}). 
\end{aligned}
\end{equation}  
The vector-valued function $y$ described above is said to be a {\it flat output} for 
\eqref{generalControlSystem}. In other words, all the trajectories of \eqref{generalControlSystem} can be expressed explicitly in terms of $y$ and finitely many of its derivatives. This is in contrast to the generic case of nonflat systems for which the representation of trajectories requires the evaluation of integrals; that is,  non--algorithmic operations. For instance, the solutions of the underdetermined ODE 
$\dot{z}=F(\ddot{y})$ can be expressed  without integrals if and only if $F''=0$ \cite{hilbert}. 

The formula
\eqref{explicitSOLformula} is said to be an {\it explicit solution} \cite{CKV24}, a {\it Monge parametrization} \cite{pomet} or a {\it flat parametrization}. Flatness is equivalent to  dynamic feedback linearizability (\cite{CKV24}, Proposition 2.1). 

Methods for the determination of flat outputs have been investigated by a number of authors. There is a significant literature on this topic; see \cite{Levine2009} and \cite{rigatosBook2015} for references up to 2015. For more recent work we mention, for instance,  \cite{gstottner}, \cite{NicRes} and the subsequent work of these authors, as  effective approaches, among others. 

In this paper we shall  take an approach which is independent of the order of the dynamic compensator by studying the  symmetries of the gUSV, namely its invariance under translations  in the plane together with a certain ``scaling'' symmetry, leading to a linearizable dynamic extension of order 6, in the general case.

 \section{Lie Symmetries of the Generic USV}\label{symmetrySect}
 
It is well known that symmetry has played a significant role in the construction of solutions to differential equations since the original work of Sophus Lie in the late 19th century.  Over the last few decades this has been a central theme in geometric mechanics \cite{holm} and geometric control \cite{vanDerSchaft}, \cite{jurdjevic}. More recently it has found application to the  study of feedback linearization in nonlinear control systems \cite{CKV24},  \cite{Klotz}, \cite{DTVa}, \cite{cascade}.  

To explain the symmetry approach to feedback linearization it is convenient to represent  
 \eqref{surfaceKinematics} as the system of differential 1-forms
\begin{equation}\label{shipPfaff}
\begin{aligned}
&\o^1=dx-(u\cos\th-v\sin\th)dt,\cr
&\o^2=dy-(u\sin\th+v\cos\th)dt,\cr
& \o^3=d\th-rdt,\cr   
&\o^4=du-\Big(\frac{1}{\g_1} vr-\b_1 u+u_1\Big)dt,\cr
&\o^5=dv+(\g_1 ur+\b_2 v)dt,\cr 
&\hskip 0 pt \o^6=dr-\left( \g_2 uv-\b_3 r+u_2\right)dt,
\end{aligned}
\end{equation}
defined on  $\B R^9$ with coordinates $t, x, y, \th, u, v, r, u_1, u_2$. The relationship between \eqref{surfaceKinematics} and \eqref{shipPfaff} is as follows: if 
$s(t)=(t,  x(t),  u(t))$ is the graph of a trajectory $(x(t), u(t))$ of \eqref{surfaceKinematics} then
$s^*\o^i=0$, $i=1,\ldots,6$, where $s^*$ denotes the {\it pullback} by $s$. 
That is, each 1-form $\o^i$ vanishes along the trajectories of \eqref{surfaceKinematics}.  Conversely, every map $s$ that satisfies $s^*\o^i=0$ for each $i$, is the graph of a trajectory of \eqref{surfaceKinematics}.
In this case  
$s:\B R\to \B R^9$, is said to be an {\it integral manifold} of
 $\bsy{\o}:=\{\o^1,\ldots,\o^6\}$.   Pictorially, $\bsy{\o}$ vanishes upon the image of  $s$ in $\B R^9$.
 
 A {\it Lie  transformation group} is a group  $\G$ of transformations 
$\t_g:M\to M$ of some manifold $M$, the elements of which are parametrized by a Lie group $G$. Each $g\in G$ defines an element $\t_g\in\G$; we then say that $G$ {\it acts on} $M$. The elements 
$\t_g\in\G$ form a group in the sense of functional composition: if $\t_{g_1}, \t_{g_2}\in \G$ then $\t_{g_1}\circ\t_{g_2}=\t_{g_1\cdot g_2}\in \G$ where $g_1\cdot g_2\in G$. 

For the model \eqref{surfaceKinematics} it is clear that there is  no preferred initial position nor initial orientation of the gUSV.  Therefore   the group of translations and rotations of the plane of the gUSV will act as its symmetries. To be clear, this consists of the transformations $\t_g$, $g=(a, b, c)$ of the configuration space $x,y,\th$ 
where  $(\bar{x}, \bar{y},\bar{\th})=\t_g(x,y,\th)$   is given by
\begin{equation}\label{symsUSV}
\begin{aligned}
&\bar{x}=x\cos a-y\cos a+b\cr
&\bar{y}=x\sin a+y\cos a+c,\ \ \ \ \ a, b, c\in \B R, \cr
&\bar{\th}=\th+a.
\end{aligned}
\end{equation}
The remaining variables $t, u, v, r, u_1, u_2$ are left untransformed. 

The transformations \eqref{symsUSV} are {\it symmetries} of \eqref{surfaceKinematics} in the sense that if $(x(t), y(t),\ldots, u_1(t), u_2(t))$ is a trajectory of  \eqref{surfaceKinematics} then its image under $\t_g$ will also be a trajectory of \eqref{surfaceKinematics}, for each $g\in G$.
 
We can express this property of mapping `trajectories to trajectories' of \eqref{surfaceKinematics}  by the {\it invariance} of $\bsy{\o}$ under $\t_g$ in the form
$$
\t_g^*\o^i\in \bsy{\o}\ \ \ \  \forall\ (a,b,c),\ 1\leq i\leq 6;
$$
that is, $\t_g^*\o^i$ is a linear combination of the elements of $\bsy{\o}$ over the smooth real-valued functions on  $\B R^9$. For if $s$ is an integral manifold of $\bsy{\o}$ then $\t_g\circ s$ is also an integral manifold since  $(\t_g\circ s)^*\bsy{\o}=s^*\t^*_g\bsy{\o}=s^*\bsy{\o}=0$.

Here, the Lie group $G$, in the discussion above,  consists of all matrices of the form
\begin{equation}
g=\left(\begin{matrix}\cos a & -\sin a & 0 & b\cr
                            \sin a &\cos a & 0 & c\cr
0 & 0 & 1 & a\cr
0 & 0& 0& 1\end{matrix}\right), \ a, b, c\in \B R,
\end{equation}
that act on the vector $\mathbf{x}=\left(\begin{matrix} x & y & \th & 1\end{matrix}\right)^{\mathrm{T}}$ by 
$$
\left(\begin{matrix}\bar{x}\cr \bar{y}\cr \bar{\th}\cr 1\end{matrix}\right)=
\left(\begin{matrix}\cos a & -\sin a & 0 & b\cr
                            \sin a &\cos a & 0 & c\cr
0 & 0 & 1 & a\cr
0 & 0& 0& 1\end{matrix}\right)
\mathbf{x}=\left(\begin{matrix} x\cos a-y\cos a+b \cr x\sin a+y\cos a+c \cr \th+a \cr 1\end{matrix}\right),
$$
where the group product $g_1\cdot g_2$ arises from the multiplication of the corresponding matrices.
This encapsulates the transformation group $\G$ given by \eqref{symsUSV} in terms of the Lie group of matrices $G$, where the group operation is matrix multiplication.

How does one discover the all the Lie symmetries of \eqref{surfaceKinematics} in the first place? The answer was fully described by Lie. In more modern terms, one seeks vector fields $X$ on $\B R^9$ that satisfy 
\begin{equation}\label{symmetryCond}
\CL_X\o^i\in\bsy{\o},\ \ 1\leq i\leq 6,
\end{equation}
where $\CL_X$ denotes the Lie derivative with respect to $X$. This leads to a set of linear homogeneous partial differential equations for the components of $X$. These are essentially solved by the Frobenius theorem. Carrying out this computation for vector fields restricted to the configuration space of the gUSV \eqref{surfaceKinematics} obtains a vector space spanned by the three vector fields:
$$
\mathfrak{g}=\{X=\P x, \ \ Y=\P y,\ \  R=\P {\th}-y\P x+x\P y\},
$$ 
that has the nonzero Lie bracket relations $[X,R]=Y,\ \  [Y,R]=-X$. The vector space
$\mathfrak{g}$ is said to be the {\it Lie algebra} of the transformation group $\G$. Then $\mathfrak{g}$  is said to constitute the {\it infinitesimal symmetries} of $\bsy{\o}$ and we say that 
$\bsy{\o}$ is {\it invariant} under $G$, or $G$-{\it invariant}.

Finally, the transformations $\t_g$ given in \eqref{symsUSV} may be constructed by composing the flow diffeomorphisms of the elements of 
$\mathfrak{g}$. In practice, however,  we usually do not need the transformations themselves, such as \eqref{symsUSV}, to study control systems -- the corresponding Lie algebras, such as $\mathfrak{g}$, are sufficient. Thus, studying the symmetries of control systems is usually computationally and conceptually straightforward. Indeed, the calculation of Lie algebras, such as $\mathfrak{g}$, are usually constructed automatically using standard procedures on a computer. In this paper we use the {\tt Maple} package {\tt DifferentialGeometry} \cite{DiffGeom}.

Our discussion above allowed us to discuss symmetries in the familiar context of rotations and translations which are obvious symmetries of the gUSV for any hydrodynamic parameter values.
However,  seeking vector fields on the entire $\B R^9$ satisfying \eqref{symmetryCond}, subject to certain technical restrictions, one finds that the generic gUSV \eqref{surfaceKinematics} has a further symmetry beyond the obvious planar Euclidean transformations, namely a scaling symmetry generated by the vector field
$$
S=x\P x+y\P y+u\P u+v\P v+u^1\P {u^1}-2\g_2uv\P {u_2},
$$
that exhibits the pivotal parameter $\g_2$.  It is {\it this} symmetry together with the translations generated by $\P x$ and $\P y$ that permits a dynamic feedback linearization that is amenable to explicit calculation. 

\section{Dynamic Feedback Linearization via Symmetry}\label{DFLsect}

In this section we derive a dynamic feedback linearization of the gUSV \eqref{surfaceKinematics} for the case $\g_2\neq 0$. The symmetries we shall study in relation to this model are generated by the Lie algebra of vector fields
$$
\bsy{\G}:=\{X, Y, S\},
$$
as discussed in Section \ref{symmetrySect}. We shall denote this Lie algebra by $\bsy{\G}$, instead of the symbol $\mathfrak{g}$ used in the previous section. The actual transformations generated by $\bsy{\G}$ are not directly relevant to our purpose and will not therefore be recorded. 

To explain our approach, recall that one possible way of constructing a dynamic feedback linearization of a control system
\begin{equation}\label{genSystem}
\dot{x}=f(t,x,u^i),
\end{equation}
which is not \sfl, is by constructing a {\it partial prolongation}, in which one or more of the inputs are differentiated some number of times:
\begin{equation}\label{ppgenSystem}
\dot{x}=f(t,x,u^i),\ \  \dot{u}^{p_i}=v^{p_i}_1, \ \dot{v}^{p_i}_1=v^{p_i}_2, \ldots, 
\dot{v}^{p_i}_{p_i-1}=v^{p_i}_{p_i}, 
\end{equation} 
in which $x, u^{p_i}$ and $v^{p_i}_{\ell_i}, 1\leq \ell_i\leq p_i-1$ are  the states of the partially prolonged system \eqref{ppgenSystem} and $v^{p_i}_{p_i}$ are its inputs, together with any of the original inputs that remain undifferentiated.

It may then happen that though \eqref{genSystem} is not \sfl, nevertheless this may hold for \eqref{ppgenSystem}. We will see that this strategy does in fact provide a dynamic feedback linearization,  for the hovercraft vessel $(\g_2=0)$ by differentiating  $u_1$ twice in \eqref{surfaceKinematics}, leading to the well-known flat outputs $[x,y]$ -- see Appendix \ref{secA2}. However, this strategy fails to lead to a \df\ linearization upon differentiating either of $u_ 1$ or $u_2$ in \eqref{surfaceKinematics}, regardless of the choice of hydrodynamic parameter values. 

Instead, to determine a \df\ linearization of the  gUSV ($\g_2\neq 0$) we will use some results from \cite{CKV24}. The overall strategy is as follows. Given an intrinsically nonlinear control system \eqref{genSystem}, a method is described in \cite{CKV24} for determining a new control system 
(denoted by $\CH_G$ in \cite{CKV24}), that is feedback equivalent to  \eqref{genSystem},  and such that a partial prolongation of 
$\CH_G$ is \sfl, even though no partial prolongation of the original system \eqref{genSystem} is \sfl. The method also  identifies which of the new inputs in $\CH_G$ are to be differentiated and  how many derivatives are required for each such input. This provides a constructive sufficient condition for \df\ linearizability.

We begin by discussing the notion of a quotient control system. Suppose a control system $\bsy{\o}$ on a manifold $M$ is $G$-invariant. That is, $G$ defines of a group of transformations with a Lie algebra 
$\bsy{\G}$ of infinitesimal generators that satisfy $\CL_X\bsy{\o}\in\bsy{\o}$ for all 
$X\in \bsy{\G}$. Provided the action of $G$ is {\it control admissible} (\cite{CKV24}, Section 3.3) then there is a control system, denoted $\bsy{\o}/G$, on the quotient of $M$ by the action of $G$, denoted $M/G$. 
The quotient space $M/G$ is defined by the equivalence relation arising from the action of $G$: two points $p,q\in M$ project to the same point in $M/G$ if and only if there is an element of $G$ that maps $p$ to $q$ under the action. For the group actions $G$ under consideration the quotient space $M/G$ is a smooth manifold and $\pi: M\to M/G$ is a smooth surjective submersion which maps each point $p\in M$ to its equivalence class $\pi(p)=[p]\in M/G$.  

To apply this construction to control systems in local coordinates we construct a local trivialization
$\bsy{u}:\pi^{-1}(\mrm U)\to \mrm U\times G$, where $\mrm U\subset M/G$ is an open set. For coordinates on $\mrm U$ we taken a collection of smooth functions $q^1, q^2, q^3, v_1, v_2$ on $M$ together with $t$
that are invariant under the action of $G$. The reason there are 6 independent invariant functions is because $\dim\B R^9-\dim G=6$, and $G$ acts {\it freely} on $\B R^9$. These functions  are easily constructed as the local first integrals of the Lie algebra
$\bsy{\G}$. It is straightforward to check that we can take
$$
t,\ \ q^1=\th,\ \  q^2=\frac{v}{u},\ \   q^3=r, \ \ v_1=\frac{u_1}{u},  \ \ v_2=u_2+\g_2 uv,
$$
for a local coordinate system on $\mrm U\subset M/G$, where $u\neq 0$. Extending these functions to a local coordinate system on $\pi^{-1}(\mrm U)$ by 
$$
g_1=x,\ \   g_2=y,  \ \ g_3=u
$$
defines the local diffeomorphism $\bsy{u}$. 

In these coordinates the map $\pi:M\to M/G$ has local expression 
$$
\pi(q^1,q^2,q^3, v_1,v_2, g_1,g_2,g_3)=(q^1,q^2,q^3, v_1,v_2)
$$
To construct the quotient control system it is convenient to pass to the smooth distribution 
\begin{equation}\label{dualShipSys}
\begin{aligned}
&\text{ann}\,\bsy{\o}=\CV=\Big\{\P t+(u\cos\th-v\sin\th)\P x+(u\sin\th+v\cos\th)\P y+r\P{\th}\cr
&+\left(\frac{1}{\g_1}vr-\b_1u+u_1\right)\P u-(\g_1ur+\b_2v)\P v+(u_2+\g_2uv-\b_3r)\P r, \P {u_1}, \P {u_2}\Big\}
\end{aligned}
\end{equation}
that is annihilated by $\bsy{\o}$. Note that is $s:\B R\to\B R^9$ is an integral manifold of $\bsy{\o}$
then it will also be an integral manifold of $\CV$ in the sense that $s_*(\P t)\in\CV_{s(t)}$. Hence,
$s$ will be an integral manifold of $\CV$ if and only if it is the graph of a trajectory of \eqref{surfaceKinematics}. Using \eqref{dualShipSys}, we get an easy representation of the quotient control system 
$\CV/G$ by the symmetry group $G$. Indeed, 
this  quotient is given by 
$$
\CV/G:=\pi_*\CV.
$$ 
A calculation shows that
$$
\CV/G=\left\{Y,\ \P {v_1},\ \P {v_2}\right\}
$$
where
$$
Y=\P t+q^3\P {q^1}-\left(\left(\b_1+\b_2-v_1\right)q^2+
\frac{\g_1^2 -(q^2)^2}{\g_1}q^3\right)\P {q^2}+\left(v_2-\b_3 q^3\right)\P {q^3}.
$$
While $\CV/G$ has a  somewhat complicated appearance, it has the crucial advantage of being \sfl, in contrast to $\CV$ itself which is not \sfl.  Indeed, the feedback transformation  $h:\mrm U\to J^{\langle 1,1\rangle}$ defined by
$$
z=q^2, z_1=v_1-\b_1+\frac{q^2q^3}{\g_1}, w=q^1, w_2=v_2-\b_3 q^3
$$
satisfies
$$
h_*(\CV/G)=\big\{\P t+z_1\P z+w_1\P w+w_2\P {w_1}, \P {z_1},\ \P {w_2}\big\}.
$$
Here the symbol $J^{\langle 1,1\rangle}$ denotes relevant (jet) space that carries the Brunosky normal form $\bsy{\b}^{\langle 1,1\rangle}=\{dz-z_1dt, dw-w_1dt, dw_1-w_2dt\}$ in which there is one chain that occurs to order 1 (the $z$-chain) and one that occurs to order 2 (the $w$-chain). For, as we shall soon see, the quotient $\CV/G$ is static feedback equivalent to 
$$
\CB_{\langle 1,1\rangle}:=\mrm {ann}\,\bsy{\b}^{\langle 1,1\rangle}=\{\P t+z_1\P z+w_1\P w+w_2\P {w_1}, \P {z_1}, \P {w_2}\}.
$$
\vskip 5 pt
Furthermore, $\bsy{u}$ maps the Lie algebra $\bsy{\G}$ to the canonical form,
$$
\bsy{u}_*\bsy{\G}=\big\{g_1\P {g_1}+g_2\P {g_2}+g_3\P {g_3}, \P {g_1},\ \P {g_2}\big\}.
$$
Next, we define the composition of feedback transformations 
$$
\vf:=(h\times \text{Id}_{G\to G})\circ \bsy{u},
$$ 
where $ \text{Id}_{G\to G}$ is the identity map on $G$. We obtain,
$$
\begin{aligned}
\CH_G:=\vf_*\CV=\{\P t+z_1\P z+w_1\P w+w_2\P {w_1}+(&\cos w-z\sin w)R_1\cr
&\hskip -80 pt+(\sin w+z\cos w)R_2-\frac{z_1+\b_2 z+\g_1 w_1}{z}R_3, \P {z_1},\ \P {w_2}\},
\end{aligned}
$$
where $R_i=g_3\P {g_i},\ 1\leq i\leq 3$. The system $\vf_*\CV$ is called the {\it contact sub-connection} 
$\CH_G$ of system \eqref{surfaceKinematics} relative to group $G$ generated by $\bsy{\G}$ (\cite{CKV24}, Theorem 4.1). As promised, $\CH_G$ is a control system on $J^{\langle 1,1\rangle}\times G$ feedback equivalent to \eqref{surfaceKinematics} that has the property of admitting a partial prolongation along $w_2$ that is \sfl. This partial prolongation is not ad-hoc and is explicitly determined by Theorem 7.6 of \cite{CKV24} using machinery we do not elaborate on here. Regardless, it can be checked that 
$$
\begin{aligned}
\text{pr}\,\CH_G=\Big\{\P t+z_1\P z+\sum_{j=0}^5w_{j+1}\P {w_j}+(\cos w-z\sin &w)R_1+(\sin w+z\cos w)R_2\cr
&-\frac{z_1+\b_2 z+\g_1 w_1}{z}R_3, \P {z_1},\ \P {w_6}\Big\},
\end{aligned}
$$
is \sfl: the signature\footnote{The notion of {\it signature} is explained in Appendix \ref{secA1}. If $\r_j$ is a nonzero entry in position $j$ of a \sfl\ control system's signature $\langle\r_1,\r_2,\ldots,\r_k\rangle$, then the corresponding Brunovsky normal form has $\r_j$ chains of order $j$.} 
of $\text{pr}\,\CH_G$ is $\langle 0,0,0,1,0,1\rangle$. In other words,
$\text{pr}\,\CH_G$ is feeback equivalent to the Brunovsky normal form on 
 the (jet) space $J^{\langle 0,0,0,1,0,1\rangle}$. That is, the Brunovsky normal form of $\mrm {pr}\,\CH_G$ has one chain of order 4 and one of order 6.
We have therefore proven
\begin{theorem}\label{gUSV-CFL-DFL}
For all parameter values $\g_1, \g_2, \b_1, \b_2, \b_3$ with $\g_2\neq 0$, the gUSV \eqref{surfaceKinematics} is \dfl, with dynamic extension of order 6.
\end{theorem}
 
We mention that one can study other symmetry quotients of the gUSV system; however, the symmetry generated by the vector field $S$ is essential to Theorem \ref{gUSV-CFL-DFL}. Indeed, Corollary 7.13 of \cite{CKV24} makes clear the important role of non-state space symmetries in generating DFLs that are not pure partial prolongation. As such, any attempt to use symmetry to construct a DFL of gUSV must somehow use the vector field $S$. 

In fact, the reader might wonder if the full four dimensional symmetry group $G'$ generated by $\bsy{\G}'=\{ X,Y,R,S\}$ must actually lead to a SFL system equivalent to the Brunovsky normal form on $J^{\langle 2 \rangle}$ (compared to the $J^{\langle 1,1 \rangle}$ for the proof of Theorem \ref{gUSV-CFL-DFL}). In fact, this turns out to be true, and would also produce a proof of Theorem \ref{gUSV-CFL-DFL}; however, care must be taken to address a curious technical phenomenon in implementing Theorem 7.6 of \cite{CKV24} to produce a minimal partial prolongation prolongation of $\text{pr}\CH_{G'}$ that is SFL. 

Finally, we remark that the question of which subgroups of a full control admissible symmetry group of a given control system admit SFL quotients remains a technical case-by-case situation. One can sometimes use only bundle rank invariants to make such deductions; see Chapter 14.9 of \cite{KlotzVassiliouQuotient} for an example classifying which control admissible subgroups of dimensions one and two yield linearizable systems. 
 
\section{Flat Parametrization of the Generic USV}\label{flatParaSect}
 
Following Theorem 1, in this section we derive a flat parametrization \eqref{explicitSOLformula} for the system \eqref{surfaceKinematics} when $\g_2\neq 0$, using the distribution $\text{pr}\CH_G$ obtained in Section \ref{DFLsect}. We also comment on the dynamic feedback linearization of  \eqref{surfaceKinematics} when $\g_2=0$ in Appendix \ref{secA2}.

We can compute  flat outputs of  \eqref{surfaceKinematics} when $\g_2\neq 0$ by using a simplified version of procedure 
{\tt contact}  \cite{Vassiliou2006b}  applied to $\text{pr}\CH_G$. However, it turns out that the \stf\ transformation to Brunovsky normal form obtained by differentiating these flat outputs is very hard to invert explicitly. As far as we can determine, there are no other choices of symmetry group that lead to a simpler \stf\ transformation for the corresponding flat outputs.  

On the other hand, examining the feedback transformation $\vf$ constructed in Section
\ref{DFLsect} shows that the new input $w_2$ in $\CH_G$ has the form
$$
w_2=v_2-\b_3q^3=u_2+\g_2uv-\b_3r,
$$
which agrees with the coefficient of $dt$ in the 1-form $\o^6$ of equation \eqref{shipPfaff}. This is the input that must be differentiated 4 times to arrive at a \sfl\ control system from  $\CH_G$. 
Thus, making the feedback transformation 
\begin{equation}\label{fbackTrans}
\bar{u}_2=u_2+\g_2uv-\b_3r
\end{equation}
in \eqref{dualShipSys} will determine a \df\ linearization. Indeed,  the four-fold partial prolongation of $\CV$ in equation \eqref{dualShipSys} along 
$\bar{u}_2$  leads to the \sfl\ control system
\begin{equation}\label{dflUSV}
\begin{aligned}
{\mrm {pr}}{\CV}=\Big\{\P t+(u&\cos\th-v\sin\th)\P x+(u\sin\th+v\cos\th)\P y+r\P{\th}\cr
&+\left(\frac{1}{\g_1}vr-\b_1u+u_1\right)\P u-(\g_1ur+\b_2v)\P v+\bar{u}_2\P r
\cr
&\hskip 40 pt +\bar{u}_{21}\P {\bar{u}_2}+\bar{u}_{22}\P {\bar{u}_{21}}+
\bar{u}_{23}\P {\bar{u}_{22}}
+\bar{u}_{24}\P {\bar{u}_{23}}, \P {u^1}, \P {\bar{u}_{24}}\Big\}.
\end{aligned}
\end{equation}
Applying {\tt contact} to ${\mrm {pr}}{\CV}$ obtains flat outputs that lead to a \stf\ linearization, whose inverse provides a flat parametrization.
 
 Using {\tt contact} we obtain  flat outputs $[y^1,y^2]$ of ${\mrm {pr}}{\CV}$ in the form
\begin{equation}\label{newOutput}
\begin{aligned}
 y^1
=rv+\frac{\g_1}{1+\g_1}\big(\begin{matrix}\, \Re && &\Im\, \end{matrix}\big) ze^{\i\th}
\left(\begin{matrix}(1+\g_1)r^2\cr \b_2r+\bar{u}_2\end{matrix}\right),
 \end{aligned}
\end{equation}
 of order 4, and
 $$
 y^2=\th,
 $$
 of order 6, where $z=x+\i y$, while $\Re$ and  $\Im$ denote the real and imaginary parts respectively, of  the complex-valued function $ze^{\i\th}$ that embodies the configuration of the non-hovercraft USV. Despite the complicated appearance of $y^1$ compared to the order 4 flat output  arising from
$\text{pr}\CH_G$, the resulting transformation that linearizes ${\mrm {pr}}{\CV}$ is amenable to explicit computation. Indeed, the Jacobian determinant of the feedback transformation determined by
$[y^1, y^2]$ is given by
$$
\Delta:=\delta\Big((1+\g_1)^2r^4-
(1+2\g_1)(\bar{u}_2)^2+(1+\g_1)r\bar{u}_{21}
+\b_2\big(\b_2r^2+(1-\g_1)r\bar{u}_2\big)\Big)^4
$$
where $\ds\delta=\frac{\g_1^4}{(\g_1+1)^5}$. The formulas for the flat parametrization are easily computed by {\tt Maple} - see Appendix \ref{secA1} and {\tt Maple} worksheet \texttt{genUSV} in Appendix \ref{secA3}. 

Observe that $y^1$ has a singularity at parameter value $\g_1=-1$.
While parameter values $\g_1\leq 0$ are unphysical,  nevertheless, the control system \eqref{surfaceKinematics} with $\g_1=-1$ is  also  \dfl\ with flat output $[y^1=y\cos\th-x\sin\th, y^2=\th]$. But the order of the dynamic extension is 4 rather than 6.

We conclude by making a brief comment on the flat parametrization obtained from the flat outputs given by \eqref{newOutput} and $y^2=\th$, referring to  the accompanying 
{\tt Maple} worksheet {\tt genUSV} for complete details. 
For notational simplicity, we express the flat parametrization in terms of the notation,
$$
z=y^1, w=y^2, z_i=\frac{d^i y^1}{dt^i}, w_j=\frac{d^j y^2}{dt^j}, \ 1\leq i\leq 4,\ \ 1\leq j\leq 6.
$$ 
Then, for instance, with $\g_1=1$, the flat parametrization for the sway velocity $v$ is given by
$$
\begin{aligned}
v=&\wt{\Delta}^{-1}\Big(\left(2w_1^5z+2w_1^3z_2-6w_1^2w_2z_1-2w_1^2w_3z+6w_1w_2^2z\right)\b_2^2\cr
&+2z\left(4w_1^4w_2+w_1^2w_4-4w_1w_2w_3+3w_2^3\right)\b_2\cr
&+\big(8zw_1^7+8z_2w_1^5-4(10w_2z_1+w_3z)w_1^4+42w_2^2zw_1^3\cr
&+4(w_3z_2-w_4z_1)w_1^2
-(6w_2^2z_2-4w_2w_3z_1-6w_2w_4z+4w_3^2z)w_1\cr
&+6w_2^2(w_2z_1-w_3z)\Big),
\end{aligned}
$$
where $\wt{\Delta}=(\b_2^2w_1^2+4w_1^4+2w_1w_3-3w_2^2)^2=
\sqrt{\phantom{(}\vf^*\Delta/\g}_{\,|_{\g_1=1}\phantom{)}}$,  $\ds \g=\g_1^4/(1+\g_1)^5$. The flat parametrization for the remaining variables is described in the {\tt Maple} worksheet {\tt genUSV}.

\backmatter

\bmhead{Supplementary information}

The calculation of flat outputs and the flat parametrization of the generic USV \eqref{surfaceKinematics} is given in the accompanying {\tt Maple} worksheet {\tt genUSV} that calls the {\tt DifferentialGeometry} package. The dynamic feedback linearization of the hovercraft is verified in the {\tt Maple} worksheet {\tt hovercraftUSV}. These worksheets are in Appendix \ref{secA3}.

\begin{appendices}

\section{Procedure {\tt Contact}}\label{secA1}

This appendix outlines the method  used for constructing the flat outputs and the \stf\ linearizations
in Section  \ref{flatParaSect}. The method is based on \cite{Vassiliou2006b}, where proofs and illustrative examples can be found. It does not assume that the system is control affine and recovers, as special cases, the known results for that class. 

We first set out some notation.
Let 
$$
\CV=\{Z:=\P t+f^j(t,x,u)\P {x^j}, \P {u_1},\ldots,\P {u^m}\}\subset TM
$$ 
\vskip 10 pt
\noindent be the vector field representation of the bracket-generating control system $\dot{x}=f(t,x,u)$ on $M$ of \textit{\textbf{derived length}} $k$. The filtration
$$
\CV\subset\CV^{(1)}\subset\CV^{(2)}\subset\cdots\subset\CV^{(k-1)}\subset\CV^{(k)}=TM
$$
is called the \textit{\textbf{derived flag}} of $\CV$, where $\CV^{(i)}=\CV^{(i-1)}+[\CV^{(i-1)}, \CV^{(i-1)}]$,
$\CV:=\CV^{(0)}$, $i\geq 1$.

\vskip 10 pt
Denote\footnote{The terms {\it bundle} and  {\it sub-bundle} (of the tangent bundle $TM$ of a manifold $M$) used here is the same as a smooth constant rank distribution over $M$.} by $\ch\CV^{(j)}$ the 
\textbf{\textit{ Cauchy bundle}} of $\CV^{(j)}$, 
$$
\ch\CV^{(j)}=\left\{X\in\CV^{(j)}~|~\big[X, \CV^{(j)}\big]\subset\CV^{(j)}\right\},\ \  j\geq 0.
$$
We assume for all $j\geq 0$, that $\CV^{(j)}$ and $\ch\CV^{(j)}$ have constant rank and refer to such sub-bundles $\CV$ as {\it totally regular}. We assume  that all sub-bundles  are totally regular. The Cauchy bundles
$\ch\CV^{(j)}$ can be shown to be  integrable for each $j$. 

Define the 
\textbf{\textit{ intersection bundles}} by
\begin{equation}\label{intersectionBundledefinition}
\ch\CV^{(i)}_{i-1}:=\CV^{(i-1)}\cap\ch\CV^{(i)},\ \ 1\leq i\leq k-1.
\end{equation}
Unlike the Cauchy bundles, the intersection bundles $\ch\CV^{(i)}_{i-1}$ are not guaranteed to be integrable in general but are integrable if $\CV$ is \sfl.

The codistribution version of Cauchy and intersection bundles respectively are denoted 
\begin{equation*}
\begin{aligned}
\Xi^{(i)}=\mathrm{ann}\,\ch\CV^{(i)},\ \ \ \ \ 
\Xi^{(i)}_{i-1}=\mathrm{ann}\,\ch\CV^{(i)}_{i-1}.
\end{aligned}
\end{equation*}

\vskip 10 pt
Let $\mathcal{V}\subset TM$ be a \sfl\ control system $\CV$ of derived length $k$ and $Z$ any vector field in $\mathcal{V}$ such that $Zt=1$. Then the \textbf{\textit{fundamental bundle}} (or \textbf{\textit{ highest order bundle}})  is given by
\begin{equation}\label{ad char}
\begin{aligned}
\Pi^{(k)}&=\big\{\Pi^0,\mathrm{ad}(Z)\Pi^0,\ldots,\mathrm{ad}(Z)^{k-1}\Pi^0 \big\},
\ \ \ \Pi^0=\ch\CV^{(1)}_0.
\end{aligned}
\end{equation}

Let $\CV\subset TM$ be a sub-bundle of derived length  $k>1$.
The \textbf{\textit{ velocity}} of $\CV$ is the ordered list of $k$ integers 
$$
\begin{aligned}
\text{\rm vel}(\CV)=\langle\D_1,\D_2,\ldots,\D_k\rangle, \ \ \ \ \ \ \text{where}\ \ \  
\D_j=\rank(\CV^{(j)})-\rank(\CV^{(j-1)}),\ 1\leq j\leq k.
\end{aligned}
$$
The \textbf{\textit{deceleration}} of $\CV$ is the ordered list of $k$ integers 
$$
\begin{aligned}
\text{\rm decel}(\CV)=\langle -\D^2_2,-\D^2_3,\ldots,-\D^2_k, \, \D_k\rangle,\ \ \ \ \ \ \text{where},\ \ \ 
\D^2_j=\D_j-\D_{j-1}.
\end{aligned}
$$
We denote the $j^{\text{th}}$ element of $\text{\rm decel}(\CV)$ by non-negative integers $\r_j$. The deceleration of
$\CV$ is also called its \textbf{\textit{signature}}. If a \sfl\ control system of derived length $k$ has signature 
$\k=\langle \r_1,\ldots,\r_k\rangle$ then it is static feedback equivalent to the unique Brunovsky normal form that has $\r_j$ chains of integrators of order $j$.

\begin{definition}[Fundamental functions, order $j<k$] \cite{Vassiliou2006b}\label{funFunction_j}
Let $\mathcal{V}\subset TM$ be a \sfl\ control system $\CV$ of derived length $k$. A set of independent first integrals 
$\big\{\phi^{\ell_j,j}\big\}_{\ell_j=1}^{\rho_j}$ of $\Xi^{(j)}_{j-1}/\Xi^{(j)}$ are called \textit{\textbf{fundamental functions of order}} $j<k$. 
\end{definition}

\begin{definition}[Fundamental functions, order $k$] \cite{Vassiliou2006b} \label{funFunction_k}Let $\mathcal{V}\subset TM$ be a \sfl\ control system $\CV$ of derived length $k$.
A set of independent first integrals 
$\big\{\phi^{\ell,k}\big\}_{\ell=1}^{\rho_k}$
of $\Pi^k$ such that $dt\w d\phi^{\ell,k}\neq 0\ \forall \ell$ are called \textit{\textbf{fundamental functions of order}} $k$.
\end{definition}

\vskip 5 pt
Let $\mathcal{V}\subset TM$ be a bracket generating \sfl\ control system of derived length $k>1$ and signature 
$\k=\langle\rho_1,\ldots,\rho_k\rangle$. 
The following steps in procedure {\bf\tt contact} provide a \stf\ linearization of $\CV$.  

\vskip 20 pt
\begin{center}-------------------------------------------------------------------------------------------\end{center}
\noindent{\large\tt \textbf{Procedure Contact}} { (for \sfl\ systems)}\ \cite{Vassiliou2006b}\label{proc B}
\begin{enumerate}
\item[] INPUT: Static feedback linearizable system $\CV=\{\P t+f^j(t,x,u)\P {x^j}, \P {u^a}\}$, with signature 
$\k=\langle \r_1,\ldots,\r_k\rangle$
\item  Fix any  vector field $Z\in \mathcal{V}$ with the property $Zt=1$. Then construct $\Pi^{k}$ as in \eqref{ad char}. 
\item For each $1\leq j\leq k-1$ such that $\rho_j\neq 0$, compute  the quotient
$\Xi^{(j)}_{j-1}/\Xi^{(j)}$.
\item For each $1\leq j\leq k-1$ such that $\rho_j\neq 0$, compute $\rho_j$ independent first integrals 
$\phi^{\ell_j,j}$ of $\Xi^{(j)}_{j-1}/\Xi^{(j)}$, $1\leq \ell_j\leq \rho_j$ ({\em fundamental functions of order $j$}).

\item Compute $\rho_k+1$ independent first integrals of the highest order bundle $\Pi^k$ and denote these
$t, \phi^{\ell_k,k}$, $1\leq \ell_k\leq \rho_k$ ({\em fundamental functions of order $k$}). 
\item For each $1\leq j\leq k$ such that $\rho_j\neq 0$, define $z^{\ell_j,j}_0=\phi^{\ell_j,j}$, $1\leq \ell_j\leq \rho_j$. The remaining contact coordinates are
\begin{equation}\label{contact coord B}
z^{\ell_j,j}_{s_j}=Zz^{\ell_j,j}_{s_j-1}=Z^{s_j}z^{\ell_j,j}_0,\,1\leq s_j\leq j,\,1\leq \ell_j\leq \rho_j. 
\end{equation}
\item[] OUTPUT: Static feedback equivalence of $\CV$ to Brunovsky normal form $\CB_\k$.
\end{enumerate}
\begin{center}-------------------------------------------------------------------------------------------\end{center}
In adapted coordinates $t,z^{\ell_j,j}_0$, together with those defined in (\ref{contact coord B}), $\mathcal{V}$ agrees with the Brunovsky normal form $\mcal{B}_\k$ of signature $\k$. 


\begin{corollary}
Let $\mathrm{pr}\,\CV$ be the \sfl\ dynamic extension  \eqref{dflUSV} of the control sytem \eqref{surfaceKinematics} represented as the distribution  $\CV$ defined in \eqref{dualShipSys}. Then the fundamental functions of
$\mrm{pr}\CV$  form a set of flat outputs for \eqref{surfaceKinematics} leading to its flat parametrization upon use of the feedback transformation \eqref{fbackTrans}.
\end{corollary}

\vskip 5 pt
\noindent By way of illustration we apply {\tt contact} to the control system $\mrm{pr}\CV$ from Section \ref{flatParaSect}, the computations for which are carried out in the {\tt Maple} worksheet 
{\tt genUSV}.

A \sfl\ dynamic extension of \eqref{surfaceKinematics} is given by the system
\eqref{dflUSV}, defined by the distribution $\mrm{pr}\CV$. As shown in equation (11) of  {\tt genUSV},  the signature of $\mrm{pr}\CV$ is $\langle 0,0,0,1,0,1\rangle$. Hence there is one fundamental function of order 4 and one of order 6. Equations (12) and (13) of {\tt genUSV} display the calculation of
$\Xi^{(4)}$ and  $\Xi^{(4)}_3$, respectively, leading to that flat output \eqref{newOutput}. Equation (18) of \texttt{genUSV} proves that \eqref{newOutput} is indeed a first integral of $\Xi^{(4)}_3$ and hence a fundamental function of order 4. Equation
(15) of  {\tt genUSV} shows that $y^2=\th$ is a fundamental function arising from the highest order bundle $\Pi^6$ and hence a flat output of order 6. The remainder of  {\tt genUSV} implements {\tt contact} to produce the flat parametrization {\tt sol} for $\mrm{pr}\CV$, just above equation (24) in {\tt genUSV},  and hence a flat parametrization for the gUSV \eqref{surfaceKinematics} via \eqref{fbackTrans}. For space reasons only the states, $x,  y, \th, v$ and $r$, of the flat parametrization, are explicitly recorded in {\tt genUSV}, as $u$ is more lengthy.

\section{Remark on the dynamic feedback linearization of the hovercraft}\label{secA2}
By way of completeness,  we briefly comment on the hovercraft vessel, which is known to be flat \cite{fliessUSV}, \cite{SiraRamirez}, in terms of our formulation.
The hovercraft is defined by \eqref{shipPfaff} with $\g_2=0$ (and hence $\g_1=1$). We denote this system  by $\bsy{\o}_0$, and  set $\CV_0:=\ker\bsy{\o}_0$. Then
$$
\begin{aligned}
&\CV_0=\Big\{Z=\P t+(u\cos\th-v\sin\th)\P x+(u\sin\th+v\cos\th)\P y+r\P {\th}\cr
&+\Big(vr-\b_1u+u_1\Big)\P u-(ur+\b_2v)\P v+(u_2-\b_3r)\P r,\ \P {u_1},\ 
\P {u_2}\Big\}.
\end{aligned}
$$
The system $\CV_0$ is given by equation (2) in the {\tt Maple} worksheet {\tt hovercraftUSV} and its prolongation along $u_1$,  $\mrm{pr}\CV_0$,
given in equation (3) of {\tt hovercraftUSV},  defined by 
$$
\mrm{pr}\CV_0=\{Z+u_{11}\P {u_1}+u_{12}\P {u_{11}}, \ \ \P {u_{12}}, \ \P {u_2}\},
$$
is \sfl,  as shown by equation (5) of {\tt hovercraftUSV}, which implements the Gardner-Shadwick-Sluis test \cite{SluisThesis}; see also \cite{KlotzVassiliouQuotient}. Finally, equation (11) of {\tt hovercraftUSV} shows that its flat outputs are $[x,y]$, which is the known result for that system.



\section{Maple files}\label{secA3}
This Appendix contains {\tt Maple} files in which the flat outputs and the flat parametrization was explicitly constructed and tested.

\end{appendices}


\bibliography{sn-bibliography}

\end{document}